\documentclass[a4paper,11pt]{my_ijamas}

\def\Sop#1{\mathrel{^*\!\mathord{#1}}}

\usepackage{colortbl}
\usepackage{graphicx}
\newcommand{\menor} {\ {\raise-.5ex\hbox{$\buildrel<\over\sim$}}\
}
\newcommand{\maior} {\ {\raise-.5ex\hbox{$\buildrel>\over\sim$}}\
}

\begin{document}


\title{Relaxed optimality conditions
for mu-differentiable functions\thanks{Accepted (18/June/2008) for
International Journal of Applied Mathematics \& Statistics
(IJAMAS).}}

\author{\textbf{Ricardo Almeida$^1$ and Delfim F. M. Torres$^2$}}

\date{$^1$Department of Mathematics\\
University of Aveiro\\ 3810-193 Aveiro, Portugal\\
ricardo.almeida@ua.pt\\[0.3cm] $^2$Department of Mathematics\\
University of Aveiro\\ 3810-193 Aveiro, Portugal\\ delfim@ua.pt}

\maketitle


\begin{abstract}
\noindent \emph{We prove some fundamental properties of
mu-differentiable functions. A new notion of local minimizer and
maximizer is introduced and several extremum conditions are
formulated using the language of nonstandard analysis. }

\medskip


\medskip


\end{abstract}


\section{Introduction}

In this work we introduce some sufficient and necessary conditions
to ensure the existence of extreme points for mu-differentiable
functions. As we will see, this type of differentiability has some
advantages when compared to others in the literature: the more
interesting one is that infinitesimal perturbations on the
function do not influence the differentiability.

The paper is organized as follows. In section \ref{NSAuniverse} we
present the usual concepts and results of Nonstandard Analysis.
The study of mu-differentiation is given in section \ref{muDif}
and in section \ref{ExtrCond} we exhibit some new (as far as we
know) sufficient and necessary conditions to guarantee the
existence of extremum points.


\section{The nonstandard universe}\label{NSAuniverse}

To prove theorems in mathematics using the $\epsilon-\delta$
definition of limit is sometimes difficult and usually not
obvious, due to the presence of three quantified, non-commutating,
expressions $\forall \epsilon \, \exists \delta \, \forall x$. In
spite of the fact that Calculus was initially formulated using
infinitesimals, in the nineteenth century mathematicians like
Augustin Cauchy, Karl Weierstrass and Richard Dedekind, working in
Mathematical Analysis, after two centuries following Isaac Newton,
usually did not mention infinitesimals. This fact lead to a 20th
century Infinitesimal Calculus where the only ``infinitesimal"
mention was in its name. However, in the past decades, things have
changed with the so called ``Nonstandard Analysis".

A number $\epsilon$ is called infinitesimal if $|\epsilon|<r$ for
all $r \in \mathbb R^+$ and $\omega=1/\epsilon$ (with $\epsilon
\not= 0$) is called an infinitely large number. In the real number
system $\mathbb R$, the only infinitesimal number is $\epsilon
=0$. However we can consider a larger system, the hyper-real
numbers ${^*\mathbb R}$, which is an ordered field that contains
$\mathbb R$ as a subfield, but also contains infinitesimals and
infinitely large numbers.

Nonstandard Analysis was invented by Abraham Robinson in the
$1960$'s, and among other things, he showed that we can embed the
ordered field of real numbers $(\mathbb R,+,\cdot,\leq)$ as an
ordered subfield of a structure $({^*\mathbb R},\Sop
+,\Sop\cdot,\Sop\leq)$ (the set of hyper-real numbers) which,
besides being a totally ordered field, contains other numbers such
as infinitesimal numbers and infinitely large numbers.

For the convenience of the reader, and in order to fix notation,
we make here a short presentation on the subject. For more about
\textit{Nonstandard Analysis} see \cite{Almeida}, \cite{Cutland},
\cite{HL}, \cite{R} and \cite{SL}.

\begin{figure}[h]
\centering
\includegraphics[width=10cm]{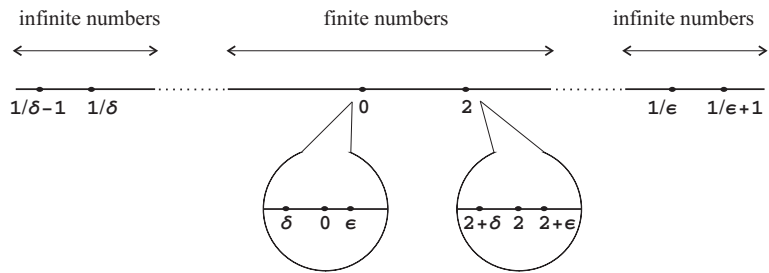}\\
  \caption{The hyper-real line.}
\end{figure}

In the following, $E$ and $F$ will denote two arbitrary (non-null)
normed spaces and ${^*E}$ and ${^*F}$ their nonstandard
extensions, respectively. These new sets contain a copy of the
primitive set
$$E \subset {^*E} \, \mbox{ and } \, F \subset {^*F}$$
but also new \textit{ideal} vectors, such as
\textit{infinitesimals}, \textit{infinite vectors}, etc (see
below).

\begin{definition} Let $x$ and $y$ be two vectors of ${^*E}$. We
say that
\begin{enumerate}
\item{$x$ is \emph{infinitesimal} if $|x|<r$ for all $r \in
    \mathbb R^+$ and we write $x\approx 0$; the set of
    infinitesimal vectors of ${^*E}$ is denoted by
    $inf({^*E})$; otherwise we write $x \not\approx 0$.}
    \item{$x$ is \emph{finite} if $|x|<r$ for some $r \in
    \mathbb R^+$ and we write $x \in fin({^*E})$.} \item{$x$
    is \emph{infinite} (or \emph{infinitely large}) if $x$ if
    not finite and we write $x\approx \infty$.} \item{$x$ and
    $y$ are \emph{infinitely close} if $x-y\approx 0$ and we
    write $x\approx y$; if not, $x \not\approx y$.} \item{$x$
    is \emph{nearstandard} ($x \in ns({^*E})$) if $x$ is
    infinitely close to some (unique) $a\in E$; in this case
    we say that $a$ is the \emph{standard part} of $x$ and we
    write $a=st(x)$.}
\end{enumerate}
\end{definition}

The rules for computing numbers in ${^*\mathbb R}$ are very simple
and they agree with our intuition. The tables below summarize
those rules. The symbols $\epsilon, a, \infty$ denote an
infinitesimal number, a finite but not infinitesimal number (what
is usually called \emph{appreciable number}) and an infinite
number, respectively (\textrm{e.g.}, $\epsilon+\epsilon=\epsilon$
means ``the sum of two infinitesimals is an infinitesimal'').

$$\begin{array}{ll}
\multicolumn{2}{l}{\,\,\curvearrowleft}\\
&\begin{array}{|c|c|c|c|} \hline
\textcolor[rgb]{0.00,0.00,0.63}{\pm}  &
\textcolor[rgb]{0.00,0.00,0.63}{\epsilon} &
\textcolor[rgb]{0.00,0.00,0.63}{a}  &
\textcolor[rgb]{0.00,0.00,0.63}{\infty} \\ \hline
\textcolor[rgb]{0.00,0.00,0.63}{\epsilon} &  \epsilon & a  &
\infty \\ \hline \textcolor[rgb]{0.00,0.00,0.63}{a}  & a  &  ?  &
\infty \\ \hline \textcolor[rgb]{0.00,0.00,0.63}{\infty} & \infty
& \infty & ?\\ \hline
\end{array}
\end{array}
\hspace{1cm}
\begin{array}{ll}
\multicolumn{2}{l}{\,\,\curvearrowleft}\\
&\begin{array}{|c|c|c|c|} \hline
\textcolor[rgb]{0.00,0.00,0.63}{\times}  &
\textcolor[rgb]{0.00,0.00,0.63}{\epsilon} &
\textcolor[rgb]{0.00,0.00,0.63}{a}  &
\textcolor[rgb]{0.00,0.00,0.63}{\infty} \\ \hline
\textcolor[rgb]{0.00,0.00,0.63}{\epsilon} &  \epsilon & \epsilon
& ? \\ \hline \textcolor[rgb]{0.00,0.00,0.63}{a}  & \epsilon  &  a
& \infty \\ \hline \textcolor[rgb]{0.00,0.00,0.63}{\infty} & ? &
\infty & \infty\\ \hline
\end{array}\end{array}
\hspace{1cm}
\begin{array}{ll}
\multicolumn{2}{l}{\,\,\curvearrowleft}\\
&\begin{array}{|c|c|c|c|} \hline
\textcolor[rgb]{0.00,0.00,0.63}{\div}  &
\textcolor[rgb]{0.00,0.00,0.63}{\epsilon} &
\textcolor[rgb]{0.00,0.00,0.63}{a}  &
\textcolor[rgb]{0.00,0.00,0.63}{\infty}\\ \hline
\textcolor[rgb]{0.00,0.00,0.63}{\epsilon} &  ? & \infty  & \infty
\\ \hline \textcolor[rgb]{0.00,0.00,0.63}{a}  & \epsilon  &  a  &
\infty \\ \hline \textcolor[rgb]{0.00,0.00,0.63}{\infty} &
\epsilon & \epsilon & ?\\ \hline
\end{array}\end{array}$$

To denote the set of infinitely large positive hyper-integers, we
use the symbol ${^*\mathbb N}_{\infty}$,
$${^*\mathbb N}_{\infty}:= \{ N \in {^*\mathbb N} \, | \, N
\approx \infty \}.$$

However, if $\omega$ and $\nu$ are two positive (resp. negative)
infinitely large numbers then $\omega + \nu$ is also a positive
(resp. negative) infinitely large number. Observe  that, if
$\epsilon \not= 0$ is an infinitesimal then
\begin{itemize}
\item $\displaystyle \frac{\epsilon^2}{\epsilon}=\epsilon$ is
    infinitesimal;\vspace{.2cm} \item $\displaystyle
    \frac{\epsilon}{\epsilon}=1$ is finite but not
    infinitesimal;\vspace{.2cm} \item $\displaystyle
    \frac{\epsilon}{\epsilon^2}=\frac{1}{\epsilon}$ is
    infinite.\vspace{.2cm}
\end{itemize}

So infinitesimals (and infinite numbers) have different orders of
magnitude. We may view the symbolic expression ``in a limit
computation, $\displaystyle\frac00$ is indeterminate" as a
short-hand for the fact that the quotient between two
infinitesimal numbers can be infinitesimal, finite but not
infinitesimal, or infinite.

It must be noted that every nearstandard vector is finite but the
opposite is false. In fact,

\begin{center}
$E$ is a finite dimensional space if and only if
$ns({^*E})=fin({^*E}).$
\end{center}

Given a vector $x \in{^*E}$, we define the \emph{monad} of $x$
has
$$\mu(x):=\{ y \in {^*E} \, | \, y-x \mbox{ is infinitesimal}
\}.$$

For example, $\mu(0)$ is simply the set of infinitesimals vectors.
For every $a \in E$, $\mu(a)=a+\mu(0)$.

\begin{theorem} Let $x$ and $y$ be two vectors in $ns({^*\mathbb
R})$. Then,
\begin{enumerate}
\item $st(x \pm y) =st(x) \pm st(y)$, $st(xy)=st(x)st(y)$ and
    $st(x/y)=st(x)/st(y)$ if $st(y) \not= 0$; \item $st(x)
    \leq st(y)$ if and only if $x \approx y$ or $x \leq y$;
    \item $st(x) \geq st(y)$ if and only if $x \approx y$ or
    $x \geq y$; \item for any $z \in \mathbb R$, $st(z)=z$;
    \item $x \approx y$ if and only if $st(x)=st(y)$.
\end{enumerate}
\end{theorem}

Clearly, some of the previous rules hold for nearstandard vectors
of ${^*E}$. For example, if $x,y \in ns({^*\mathbb E})$, it is
true that $st(x \pm y)=st(x) \pm st(y)$.

Let $U$ be a nonempty subset of $E$. The set of \emph{nearstandard
vectors} of $U$, denoted by $ns({^*U})$, is given by
$$ns({^*U}):= \{ x \in {^*U} \, | \, x \in ns({^*E}) \, \mbox{ and
} \, st(x) \in U \}.$$

\begin{definition} Given an internal function $f:{^*U}\to{^*F}$,
we say that f is \emph{S-continuous}
at $a\in {^*U}$ if  for all $x\in{^*U}$, if $x \approx a$ then
$f(x) \approx f(a)$. If $f$ is S-continuous at all $a\in U$, we
say that $f$ is \emph{S-continuous}. If it still holds for all
$a\in{^*U}$, then $f$ is said to be \emph{SU-continuous}.
\end{definition}

There exists a relation between S-continuity and (classical)
continuity. It can be proven the following.

\begin{theorem} A standard function $f:U\to F$ is continuous
(resp. uniformly continuous) if and only if its nonstandard
extension ${^*f}:{^*U}\to{^*F}$ is S-continuous (resp.
SU-continuous).
\end{theorem}

This is a very nice characterization of continuity. For example,
let us prove that $f(x)=x^2, \, x \in \mathbb R$, is continuous
but not uniformly continuous. Fix $a \in \mathbb R$. It is enough
to see that, given any infinitesimal $\epsilon$,
$$f(a+\epsilon)=a^2+2\epsilon a+\epsilon^2\approx a^2=f(a);$$
but given (any) infinite number $\omega$,
$$f \left( \omega+\frac{1}{\omega} \right)
=\omega^2+2+\frac{1}{\omega^2} \not\approx \omega^2=f(\omega).$$
If we consider, however, $f(x)=x^2, \, x \in [-1,1]$, then $f$ is
uniformly continuous. In fact, if $x,y \in {^*[-1,1]}$ with $x
\approx y$,
$$f(x)-f(y)=x^2-y^2=(x-y)(x+y)\approx 0$$
since $x-y \approx 0$ and $x+y$ is finite.

\section{The mu-differentiability}\label{muDif}

We now present the basic properties of a recent notion of
differentiation, called \emph{mu-differentiabi-} \emph{lity}
(\cite{Almeida} and \cite{AlmeidaNeves}). The important about this
type of derivative is that, under some assumptions, if $g$ is a
standard $C^1$ function, $f$ is an internal function, and $g$ is
somewhat infinitely close to $f$, then $f$ is mu-differentiable,
and vice-versa. From this we see that for standard functions,
mu-differentiability is equivalent to Fr\'{e}chet differentiation.
Therefore, the novelties appear when we work with internal (but
not standard) functions. To begin with, we will define what is the
\emph{standard part of a function} $f$. Let $f:{^*U}\to {^*F}$ be
an internal function such that $f(ns({^*U})) \subseteq ns({^*F})$.
Then, we can define a new (standard) function, which we denote by
$st(f)$,
$$\begin{array}{cccl}
st(f): & U & \to & F\\
      & x &\mapsto & st(f(x)).
\end{array}$$

Observe that, by definition, por all $x \in U$,
$st(f)(x)=st(f(x))$. For example, let
\begin{equation}\label{example1}\begin{array}{cccl}
f: & {^*\mathbb R} & \to & {^*\mathbb R}\\
      & x &\mapsto & \left\{
\begin{array}{lll}
x^2 & \mbox{ if } & x \not= 0\\ \epsilon & \mbox{ if } & x=0\\
\end{array}
      \right.
\end{array}\end{equation}

where $\epsilon$ is any nonzero infinitesimal. Then $st(f)$ is
simply  the quadratic function $x\mapsto x^2$.  Notice that the
nonstandard extension of $st(f)$ can be distinct of $f$. In this
example, since $st(f)(x)=x^2 \, , x \in \mathbb R$, then
${^*(st(f))}(x)=x^2 \, , x \in {^*\mathbb R}$.

In \cite{Schlesinger} it is presented a new kind of
differentiation:

\begin{definition} Let $U\subseteq E$ be an open set and
$f:{^*U}\to{^*F}$ be an internal function such that
$f(ns({^*U}))\subseteq ns({^*F})$. We say that $f$ is
\emph{m-differentiable} at $a \in U$ if there exists a positive
infinitesimal $\delta_a$ and an internal finite\footnote{By finite
we mean $Df_a(fin({^*E}))\subseteq fin({^*F})$.} linear operator
$Df_a \in {^*L(E,F)}$ such that
$$\forall x \in {^*U} \hspace{.5cm} \delta_a< |x-a|\approx 0
\Rightarrow f(x)-f(a)=Df_a(x-a)+|x-a|\eta,$$
for some $\eta \approx 0$ (which depends on $x$). If $f$ is
m-differentiable at all $a\in U$, we say that $f$ is
m-differentiable.
\end{definition}

There is an important criterium to test m-differentiability,
comparing with standard functions:

\begin{theorem} \cite{Schlesinger} Let $E$ and $F$ be standard
finite dimensional normed spaces, $K$ a standard compact subset of
$E$ and $f:{^*K} \rightarrow {^*F}$ an internal function. Then,
the two following statements are equivalent:
\begin{enumerate}
\item{$f$ is S-continuous and m-differentiable;} \item{There
    exists a differentiable standard function $g:K \rightarrow
    F$ such that
$$\displaystyle \sup_{x \in {^*K}}|f(x)-g(x)|\approx 0.$$}
\end{enumerate}
\end{theorem}

Basically, this result asserts that if $g$ is a standard
differentiable function and $f$ is an internal function infinitely
close to $g$, then $f$ is m-differentiable. For example, if $f$ is
the function defined above (see (\ref{example1})), and if we
define $g(x)=x^2 \, ,x \in \mathbb R$, since $f(x)\approx g(x)$
for all $x \in {^*\mathbb R}$ then $f$ is m-differentiable.

In \cite{AlmeidaNeves} this notion is extended by introducing the
concept of m-uniformly differentiability (shortly
mu-differentiable). In the following we present the main results
of that paper.

\begin{definition}\label{defmu} Let $U\subseteq E$ be an open set
and $f:{^*U} \rightarrow {^*F}$ be an internal function satisfying
$f(ns({^*U})) \subseteq ns({^*F})$. We say that $f$ is
\emph{mu-differentiable} if there exists an internal function from
$^*\!U $ into ${^*L(E,F)}$, $x\mapsto Df_x$ such that
\begin{enumerate}
\item when $x$ is near-standard in $^*U$, $Df_x$ is a finite
    map. \item for each $a \in U$, there exists a positive
    infinitesimal $\delta_a$ for which, when $x,y\approx a$,
    some infinitesimal  vector $\eta$ verifies
$$|x-y|>\delta_a \Rightarrow f(x)-f(y)=Df_x(x-y)+|x-y|\eta.$$
\end{enumerate}
\end{definition}

Actually one encompassing $\delta$ may be taken in Definition
\ref{defmu}, i.e., the following holds.

\begin{theorem}\label{tdefmu}
Let $f:{^*U} \rightarrow {^*F}$ be an internal function; $f$ is
mu-differentiable  if and only if  all the following conditions
are verified
\begin{enumerate}
\item $f(ns({^*U})) \subseteq ns({^*F})$. \item There exist an
    internal function from  $^*\!U$ into $^*L(E,F)$, $x\mapsto
    Df_x$ and a positive infinitesimal $\delta$ such that
\begin{enumerate}
\item when $x$ is near-standard in $^*U$, $Df_x$ is a
    finite map. \item when $x$  and $y$ are near-standard
    in $^*U$,  some  infinitesimal  vector $\eta$
    verifies
$$
|x-y|>\delta \Rightarrow f(x)-f(y)=Df_x(x-y)+|x-y|\eta \,
.
$$
\end{enumerate}
\end{enumerate}
\end{theorem}

The term m-uniform differentiability is justified by the following
result:

\begin{theorem}\label{theorem1} Let $f:{^*U} \rightarrow {^*F}$ be
an internal function. Then:
\begin{enumerate}
\item{If $F$ is a finite dimensional space and $f$ is a
    mu-differentiable function, then $st(f):U \rightarrow F$
    is a $C^1$ function and $D st(f)_a=st( Df_a)$ for $a \in
    U$. Furthermore, if $E$ is also  finite dimensional then
    $$\forall a \in U\, \exists \eta _0 \approx 0 \, \forall x
    \approx a \hspace{.5cm} |f(x)-st(f)(x)| \leq \eta_0.$$}
    \item{If there exists a $C^1$ standard function $g:U
    \rightarrow F$ with
$$\forall a \in U\, \exists \eta _0  \approx 0 \, \forall x
\approx a \hspace{.5cm} |f(x)-g(x)| \leq \eta_0,$$ then $f$ is
mu-differentiable. Moreover, $g=st(f)$.}
\end{enumerate}
\end{theorem}

From the previous theorem, one can prove that for standard
functions $f$,
\begin{center} $f$ is of class $C^1$ if and only if $f$ is
mu-differentiable.
\end{center}

It is clear that, if $f$ and $g$ are two mu-differentiable
functions and $k \in fin({^*\mathbb{R}})$, then $f+g$ and $kf$ are
also mu-differentiable and $D(f+g)_x=Df_x+Dg_x$ and $D(kf)_x=k \,
Df_x$.

There is some form of S-continuity for the function $f$ and for
its derivative map $Df_{(\cdot)}$:

\begin{theorem} If $f:{^*U} \rightarrow {^*F}$ is a
mu-differentiable function, then
$$\forall x,y \in ns({^*U}) \hspace{0.5cm} x \approx y \Rightarrow
f(x) \approx f(y)$$
and for all unit vector $d \in {^*E}$,
$$\forall x,y \in ns({^*U}) \hspace{0.5cm} x \approx y \Rightarrow
Df_x(d) \approx Df_y(d).$$
\end{theorem}

\begin{theorem}[Chain Rule] Let $g$ and $f$ be two
m-differentiable functions at $a$ and $g(a)$, respectively, where
$a$ and $g(a)$ are two standard vectors. In addition, if $Dg_a$ is
invertible and $\|(Dg_a)^{-1}\|$
is finite, then $f \circ g$ is m-differentiable at $a$ and
$D(f\circ g) _a=Df_{g(a)} \circ Dg_a$.
\end{theorem}

Note that since a mu-differentiable function is m-differentiable,
this result follows for mu-differen-tiability.

Next we present a Mean Value Theorem for mu-differentiable
functions. In opposite to standard functions, for an internal
function the derivative is not unique. For example, if
$f:{^*\mathbb R}\to {^*\mathbb R}$ is an internal function,
$a\in\mathbb R$ and $f'(a)$ is one derivative of $f$ at $a$ then
$f'(a)+\epsilon$ with $\epsilon\approx0$ is also a possible
derivative:
$$\delta\approx 0 \Rightarrow  f'(a) \approx
\frac{f(a+\delta)-f(a)}{\delta} \approx f'(a)+\epsilon.$$
This justifies the presence of the infinitesimal term $|x-y|\eta$
in the following result.

\begin{theorem}[Mean Value Theorem]Let $U$ be a standard open
convex subset of $E$ and $f:{^*U}  \rightarrow {^*\mathbb{R}}$ an
internal mu-differentiable function. Take $\delta$ as
given by Theorem \ref{tdefmu}. Then, for all $x,y \in ns({^*U})$
with $|x-y|>\delta $,
$$\exists c \in {[x,y]} \hspace{0.3cm} f(x)-f(y)=Df_c(x-y)+|x-y|
\eta$$
for some $\eta \approx 0.$ More generally, let $f:{^*U}
\rightarrow {^*F}$ be an internal  mu-differentiable function and
$\delta$ as given by Theorem \ref{tdefmu}. Then, for all $x,y \in
ns({^*U})$ with $|x-y|>\delta $,
$$\exists c \in {[x,y]}  \hspace{0.3cm}  |f(x)-f(y)| \leq
|Df_c(x-y)|+|x-y|\eta$$
for some $\eta \approx 0$.
\end{theorem}

A full version of an Inverse Mapping Theorem for mu-differentiable
functions is not expected. The argument is simple: the standard
function $g(x)=x \, , x \in \mathbb R$, is of class $C^1$ and
invertible. We have proved that any internal function $f$
infinitely close to $g$ is mu-differentiable, so the one-to-one
condition may easily fail. Nevertheless, we have some form of
injectivity:

\begin{theorem}[Inverse Mapping Theorem] Let $f:{^*U} \rightarrow
{^*F}$ be an internal mu-differentiable function. Assume that, for
a certain $a \in U$, $Df_a$ is invertible and $\|(Df_a)^{-1}\|$ is
finite. Then, there exists a standard neighborhood ${^*V}$ of $a$
such that $f$ is injective on the standard elements of ${^*V}$,
\textit{i.e.},
$$\forall x,y \in V \hspace{.2cm} x \not= y \Rightarrow f(x) \not=
f(y).$$
\end{theorem}

Let $f$ be a mu-differentiable function and
$$Df_{(\cdot)}:{^*U}\rightarrow {^*L(E,F)}$$
its derivative map. Since $L(E,F)$ is a standard normed space, it
makes sense to define higher-order derivatives. We say that $f$ is
twice mu-differentiable provided $f$ and $Df_{(\cdot)}$ are both
mu-differentiable.

Recursively, $f$ is $k$-times mu-differentiable ($k \in \mathbb
N$) provided $f$, $Df_{(\cdot)}$, ...,$D^{k-1}f_{(\cdot)}$ are all
mu-differentiable.

\begin{theorem}\label{theorem2} Let $f:{^*U} \rightarrow {^*F}$ be
an internal function. Then:
\begin{enumerate}
\item{If $F$ is a finite dimensional space and $f$ is
    $k$-times mu-differentiable, then $st(f):U \rightarrow F$
    is a $C^k$ function and for each $a \in U$, $D^j
    st(f)_a=st( D^jf_a)$ for $j=1,2, \ldots , k$. Furthermore,
    if $E$ is also finite dimensional,
$$\forall a \in U \, \exists \eta _0 \approx 0  \, \forall x
\approx a \hspace{.3cm} |f(x)-st(f)(x)| \leq \eta_0$$
and
$$\forall j \in \{ 1,2, \ldots, k-1 \} \hspace{.1cm}  \forall
a \in U \hspace{.1cm}  \exists \eta _j \approx 0 \hspace{.1cm}
\forall x \approx a \hspace{.5cm} |D^jf_x-D^j st(f)_x| \leq
\eta_j.$$}
\item{If there exists a $C^k$ standard function $g:U
    \rightarrow F$ with
$$\forall a \in U \, \exists \eta _0 \approx 0 \, \forall x
\approx a \hspace{.3cm} |f(x)-g(x)| \leq \eta_0$$
and
$$\forall j \in \{ 1,2, \ldots, k-1 \}  \, \forall a \in U  \,
\exists \eta _j \approx 0 \, \forall x \approx a \hspace{.3cm}
|D^jf_x-D^jg_x| \leq \eta_j$$
then $f$ is $k$-times mu-differentiable. Moreover,
$g=st(f)$.}
\end{enumerate}
\end{theorem}

\begin{theorem}[Taylor's Theorem]\label{TaylorTheor}  Let $E$ and
$F$ be two standard  finite dimensional spaces, $U\subset E$ a
standard   open set and $f:{^*U} \rightarrow {^*F}$ an internal
function $k$-times mu-differentiable, for some $k \in \mathbb{N}$.
Then,
\begin{enumerate}
\item{for every $x \in ns({^*U})$, there exists $\epsilon
    \approx 0$ such that, whenever $y \in {^*U}$ with
    $\epsilon < |y-x| \approx 0$, there exists $\eta \approx
    0$ satisfying
$$f(y)=f(x)+Df_x(y-x)+\frac{1}{2!}D^{2}f_x(y-x)^{(2)}+...+\frac{1}{k!}D^{k}f_x(y-x)^{(k)}+|y-x|^k
\eta.$$}
\item{for every $x \in ns({^*U})$, there exists $\epsilon
    \approx 0$ such that, whenever $y \in {^*U}$ with
    $\epsilon < |y-x| \approx 0$, there exists $\eta \approx
    0$ satisfying
$$f(y)=st(f)(x)+Dst(f)_x(y-x)+\frac{1}{2!}D^{2}st(f)_x(y-x)^{(2)}+...$$
$$+\frac{1}{k!}D^{k}st(f)_x(y-x)^{(k)}+|y-x|^k \eta.$$}
\end{enumerate}
\end{theorem}

\section{Main results: extremum conditions for mu-differentiable
functions}\label{ExtrCond}

In the following, $f:{^*\mathbb R}\to{^*\mathbb R}$ denotes an
internal function. The goal is to present a notion of minimizer
and maximizer, and to infer the necessary and sufficient
conditions of optimality. Obviously, the (usual) definition of
minimizer is not a good one for our study. Let
$g(x)=x^2,x\in\mathbb R$. By Theorem~\ref{theorem1}, any internal
function $f:{^*\mathbb R}\to{^*\mathbb R}$ infinitely close to $g$
is mu-differentiable and $(st(f))^{(k)}(a)=st(f^{(k)})(a)$, $a\in
\mathbb R$. Consequently, if $\epsilon$ is a positive
infinitesimal, the functions
$$\begin{array}{lcccr}
f_1(x)=x^2 \, , & \hspace{.5cm} f_2(x)=  \left\{
\begin{array}{lll}
x^2 \, & \mbox{ if } & x \not=0\\ \epsilon \, &\mbox{ if } & x
=0\\
\end{array}\right. \, ,& \hspace{.5cm}
f_3(x)=  \left\{
\begin{array}{lll}
x^2 \, &\mbox{ if } & x \not=0\\ -\epsilon \, &\mbox{ if } & x
=0\\
\end{array}  \right.
\end{array}$$
have the same derivatives of all orders (or we can choose as
such). Consequently, the definition of minimizer must take into
account this fact. In the following, given $x,y\in{^*\mathbb R}$,
$x \maior y$ (resp. $x\menor y$) will mean $x \geq y$ or $x\approx
y$ (resp. $x \leq y$ or $x\approx y$). Moreover, $x \gg y$ (resp.
$x \ll y$) is an abbreviation for $x >y$ and $x \not\approx y$
(resp. $x<y$ and $x \not\approx y$).

\begin{definition} Let $a\in\mathbb R$ be a real. We say that $a$
is a \emph{local m-minimizer} of $f$ if there exists a positive
$r\in\mathbb R$ such that
$$\forall x \in {^*]a-r,a+r[} \hspace{.5cm} f(x) \maior f(a).$$
\end{definition}

\begin{definition} Let $a\in\mathbb R$ be a real. We say that $a$
is a \emph
{local m-maximizer} of $f$ if there exists a positive $r\in\mathbb
R$ such that
$$\forall x \in {^*]a-r,a+r[} \hspace{.5cm} f(x) \menor f(a).$$
\end{definition}

Replacing $f$ by $-f$, all results proved henceforth about
m-minimums have then an equivalent for m-maximums. Without loss of
generality, from now on we will simply say that $a$ is a
m-minimizer.

The next theorem establish a relation between the standard and the
nonstandard universes. As we will see, to prove theorems we
transfer some properties to the standard universe, apply the
well-known results about standard functions and then go back to
the nonstandard universe.

\begin{lemma}\label{ponte} Let $f$ be a mu-differentiable
function. Then,
$$a \mbox{ is a m-minimizer of } f \mbox{ if and only if } a
\mbox{ is a minimizer of } st(f).$$
\end{lemma}

\begin{proof} First suppose that $a$ is a m-minimizer. Then,
$$\forall x \in {^*]a-r,a+r[} \hspace{.5cm} f(x) \maior f(a)$$
and so $st(f(x)) \geq st(f(a))$. In particular,
$$\forall x \in ]a-r,a+r[ \hspace{.5cm} st(f(x)) \geq st(f(a)),$$
\textit{i.e.}, $st(f)(x) \geq st(f)(a)$ and we proved that $a$ is
a minimizer of $st(f)$.

To prove the converse, assume
$$\forall x \in ]a-r,a+r[ \hspace{.5cm} st(f)(x) \geq st(f)(a).$$
By the Transfer Principle (see \textrm{e.g.} \cite{HL}), it holds
$$\forall x \in {^*]a-r,a+r[} \hspace{.5cm} st(f)(x) \geq
st(f)(a).$$
Since
\begin{itemize}
\item{$st(f)(a)=st(f(a)) \approx f(a)$ since $a$ is standard;}
    \item{$st(f)(x)\approx st(f) (st(x)) = st(f(st(x)))
    \approx f(st(x))\approx f(x)$ since $f$ and $st(f)$ are
    S-continuous;}
\end{itemize}
it follows
$$\forall x \in {^*]a-r,a+r[} \hspace{.5cm} f(x) \maior f(a).$$
\end{proof}

Recall Theorem \ref{theorem2}: if $f:{^*\mathbb R}\to{^*\mathbb
R}$ is $k$-times mu-differentiable, then $st(f)$ is of class $C^k$
and $(st(f))^{(k)}(a)=st(f^{(k)})(a)$, for every $a \in \mathbb
R$. Moreover, by definition of standard part of a function, if
$a\in \mathbb R$ then $st(f^{(k)})(a)=st(f^{(k)}(a))$; so if
$f^{(k)}(a)\approx L$, for some $L \in \mathbb R$, then
$st(f^{(k)})(a)=L$.

\begin{theorem}[Necessary condition for m-minimum] If $a$ is a
m-minimizer of $f$, then $f'(a)\approx 0$.
\end{theorem}

\begin{proof} If $a$ is a m-minimizer of $f$, then $a$ is a
minimizer of $st(f)$. Consequently $(st(f))'(a)=0$. By
Theorem~\ref{theorem1}, $st(f')(a)=0$ and so $f'(a)\approx 0$.
\end{proof}

\begin{theorem}[Sufficient condition for m-minimum] If $f$ is
twice mu-differentiable, $f'(a)\approx 0$, and $f''(a) \gg 0$,
then $a$ is a m-minimizer of $f$.
\end{theorem}

\begin{proof} Since $st(f'(a))=0$ and $st(f''(a))>0$, it follows
that $(st(f))'(a)=0$ and $(st(f))''(a)>0$. Then $a$ is a minimizer
of $st(f)$ and so a m-minimizer of $f$.
\end{proof}

We remark that we can prove the previous theorem using only the
Taylor's Theorem (Theorem~\ref{TaylorTheor}), avoiding the usage
of Lemma \ref{ponte}. By Taylor's Theorem, there exists
$\epsilon\approx0$ such that for all $x\approx a$, if
$|x-a|>\epsilon$ then for some infinitesimal $\eta$,
$$f(x)-f(a)=f'(a)(x-a)+f''(a)\frac{(x-a)^2}{2}+(x-a)^2\eta=(x-a)^2\left(\frac{f'(a)}{x-a}
+\frac{f''(a)}{2}+\eta \right).$$
We may also assume that $\epsilon>\sqrt{|f'(a)|}$. Then,
$$\left| \frac{f'(a)}{x-a}
\right|<\frac{|f'(a)|}{\epsilon}<\frac{|f'(a)|}{\sqrt{|f'(a)|}}\approx
0$$
and so
$$\frac{f'(a)}{x-a} +\frac{f''(a)}{2}+\eta \approx
\frac{f''(a)}{2}>0.$$
We proved that
$$\forall x\in {^*\mathbb R} \hspace{.5cm} [x\approx a \wedge
|x-a|>\epsilon] \Rightarrow f(x)>f(a).$$
Define the (internal) set $C$ as being
$$C:= \{ \theta \in {^*\mathbb R^+} \, | \, \theta \leq \epsilon
\vee (\forall \xi \in ]\epsilon,\theta[ \hspace{.5cm}
f(a\pm\xi)>f(a) ) \}.$$
Since $C$ contains all positive infinitesimal numbers, it also
contains a real $r$ (\textit{Cauchy's Principal}, see \cite{HL}).
Since $r > \epsilon$, it follows that
$$\forall \xi \in ]\epsilon,r[ \hspace{.5cm} f(a\pm\xi)>f(a).$$
For $ \xi \in [-\epsilon,\epsilon]$, $f(a+\xi)\approx f(a)$. In
conclusion,
$$\forall \xi \in {^*]-r,r[} \hspace{.5cm} f(a+\xi)\maior f(a).$$

\begin{theorem}[Higher-order necessary condition for m-minimum]
Let $f$ be a function $k$-times mu-differentiable. If $a$ is a
m-minimizer of $f$ and
$$f'(a)\approx f''(a) \approx \ldots \approx f^{(k-1)}(a)\approx
0$$
then
\begin{enumerate}
\item $f^{(k)}(a)\approx0$ if $k$ is odd; \item
    $f^{(k)}(a)\maior 0$ if $k$ is even.
\end{enumerate}
\end{theorem}

\begin{proof} If $a$ is a m-minimizer of $f$ then $a$ is a
minimizer of $st(f)$. On the other hand,
$$(st(f))'(a)=(st(f))''(a)= \ldots =(st(f))^{(k-1)}(a) =0.$$
Therefore, (see \textit{e.g.} \cite{fenske})
\begin{enumerate}
\item $(st(f))^{(k)}(a) =0 \Leftrightarrow st(f^{(k)}(a))=0
    \Leftrightarrow f^{(k)}(a) \approx 0$ if $k$ is odd; \item
    $(st(f))^{(k)}(a) \geq 0 \Leftrightarrow st(f^{(k)}(a))
    \geq 0 \Leftrightarrow f^{(k)}(a) \maior 0$ if $k$ is
    even.
\end{enumerate}
\end{proof}

\begin{theorem}[Higher-order sufficient condition for m-minimum]
If $f$ is $k$-times mu-differentiable,
$$f'(a)\approx f''(a) \approx \ldots \approx f^{(k-1)}(a)\approx
0$$ and $f^{(k)}(a)\not\approx 0$, then
\begin{enumerate}
\item if $k$ is odd then $a$ is not a m-minimizer of $f$.
    \item if $k$ is even and $f^{(k)}(a) \gg 0$, then $a$ is a
    m-minimizer of $f$.
\end{enumerate}
\end{theorem}

\begin{proof} Since
$$(st(f))'(a)=(st(f))''(a)= \ldots =(st(f))^{(k-1)}(a) =0 \mbox{
and } (st(f))^{(k)}(a)\not=0,$$
it follows that (see \textit{e.g.} \cite{fenske})
\begin{enumerate}
\item{for $k$ odd, $a$ is not a minimizer of $st(f)$ and so
    not a m-minimizer of $f$.} \item{for $k$ even, since
    $(st(f))^{(k)}(a)>0$, $a$ is a minimizer of $st(f)$  and
    so a m-minimizer of $f$.}
\end{enumerate}
\end{proof}

For example, let
$f(x)=\frac{1}{7}x^7-\frac{1}{2}x^6+\frac{2}{5}x^5+\epsilon x+
\delta x^2$, where $\epsilon$ and $\delta$ are two non-zero
infinitesimals. Then, $f'(x)=x^6-3x^5+2x^4+\epsilon+2\delta x$ and
so
$$f'(0)=\epsilon\approx0 \, , \, f'(1)=\epsilon+2\delta\approx0 \,
, \,f'(2)=\epsilon+4\delta\approx0.$$
The second derivative of $f$ is $f''(x)=6x^5-15x^4+8x^3+2\delta$.
Thus,
$$f''(0)=2\delta\approx0 \, , \, f''(1)=-1+2\delta \ll 0 \, ,
\,f''(2)=16+2\delta \gg 0.$$
Furthermore, since $f'''(0)=f^{(4)}(0)=0$ and $f^{(5)}(0) =48 \gg
0$, we conclude that
$$x=1 \mbox{ is a m-maximizer of }f \, , \, x=2 \mbox{ is a
m-minimizer of }f$$
$$\mbox{and } x=0 \mbox{ is neither a m-maximizer nor a
m-minimizer of }f.$$

\subsection*{Functions with several variables}

From now on we will work with internal functions with several
variables, $f:{^*\mathbb R}^n \to {^*\mathbb R}$. Suppose that $f$
is mu-differentiable, \textit{i.e.}, given $a\in \mathbb R^n$,
there exists some $\delta\approx 0$ such that, for all $x\approx
a$ and $\epsilon\approx0$, if $|\epsilon|>\delta$ then
\begin{equation}\label{eq1}
f(x+\epsilon)-f(x)=Df_x(\epsilon)+|\epsilon|\eta
\end{equation}
for some $\eta\approx 0$. Let $x:=(x_1,\ldots,x_n)$,
$\epsilon:=(\epsilon_1,\ldots,\epsilon_n)$ and
$Df_x:=(f'_1,\ldots,f'_n)$. Rewriting equation (\ref{eq1}),
$$f(x_1+\epsilon_1,\ldots,x_n+\epsilon_n)-f(x_1,\ldots,x_n)=f'_1\epsilon_1+\ldots+
f'_n\epsilon_n +\sqrt{\epsilon_1^2+\ldots + \epsilon_n^2} \;
\eta.$$
Consequently, for each $i \in \{1,\ldots, n \}$, if we fix
$\epsilon_i\not=0$ and $\epsilon_k=0$ for $k \in \{1,\ldots,
i-1,i+1,\ldots,n \}$,
$$f(x_1,\ldots, x_{i-1}, x_i+\epsilon_i, x_{i+1},
\ldots,x_n)-f(x_1,\ldots,x_n)=f'_i\epsilon_i+ |\epsilon_i|\eta$$
if $|\epsilon_i|>\delta$, \textit{i.e.},
$$f'_i\approx \frac{f(x_1,\ldots, x_{i-1}, x_i+\epsilon_i,
x_{i+1}, \ldots,x_n)-f(x_1,\ldots,x_n)}{\epsilon_i}.$$
Let us denote $f'_i$ by $\displaystyle \left. \frac{\partial
f}{\partial x_i} \right|_{x}$ and we call them the \textit{partial
derivatives} of $f$.

By Theorem \ref{theorem1}, $st(Df_a)=Dst(f)_a$ whenever $a$ is
standard. So
$$st\left(\left. \frac{\partial f}{\partial x_1}
\right|_{a},\ldots, \left. \frac{\partial f}{\partial x_n}
\right|_{a}\right)=\left(\left. \frac{\partial st(f)}{\partial
x_1} \right|_{a},\ldots,\left. \frac{\partial st(f)}{\partial x_n}
\right|_{a}\right)$$
\begin{equation}\label{eq3}\Leftrightarrow st \left( \left.
\frac{\partial f}{\partial x_i} \right|_{a} \right) = \left.
\frac{\partial st(f)}{\partial x_i} \right|_{a} \, \mbox{ for } \,
i \in \{1,\ldots,n \}.\end{equation}

For example, let
\begin{equation}\label{eq2}\begin{array}{lcll}
f: & {^*\mathbb R^2} & \to & {^*\mathbb R}\\
   & (x_1,x_2)      & \mapsto & \frac{\sin(\epsilon
x_1)}{\epsilon}+x_2\\
\end{array}
\end{equation}
where $\epsilon$ is a positive infinitesimal number. Let us prove
that $f$ is mu-differentiable. Denote $f$ as the sum of $f_1$ and
$f_2$. Clearly, $f_2(x_1,x_2):=x_2$ is mu-differentiable, so we
only will prove that $f_1(x_1,x_2):= \frac{\sin(\epsilon
x_1)}{\epsilon}$ is also mu-differentiable. First, observe that
for every $(x_1,x_2)\in ns( {^*\mathbb R^2})$,
$$\frac{\sin(\epsilon x_1)}{\epsilon}= \frac{\sin(\epsilon
x_1)}{\epsilon x_1} x_1 \approx x_1\in ns({^*\mathbb R}).$$
Thus $f(ns( {^*\mathbb R^2}))\subseteq ns({^*\mathbb R})$. Let
$$\begin{array}{lcll}
g: & \mathbb R^2 & \to & \mathbb R\\
   & (x_1,x_2)      & \mapsto & x_1.\\
\end{array}$$
Then $g$ is a (standard) $C^1$ function and, given $a \in \mathbb
R^2$, let $\eta:=\epsilon$. Then, for all $(x_1,x_2)\approx a$,
$$|f_1(x_1,x_2)-g(x_1,x_2)|=\left| \frac{\sin(\epsilon
x_1)}{\epsilon} -x_1 \right|=\left| \frac{\epsilon x_1 -
\frac{(\epsilon x_1)^3}{6}+\epsilon^3 \xi}{\epsilon} -x_1
\right|$$
$$=\left|- \frac{\epsilon^2 x_1^3}{6}+\epsilon^2
\xi\right|=\epsilon \left| - \frac{\epsilon x_1^3}{6}+\epsilon
\xi\right| < \eta$$
for some $\xi\approx 0$ (note that, since $x\mapsto \sin x$ is a
$C^1$ standard function, by Taylor's Theorem,
$\sin(x)=x-\frac{x^3}{6}+x^3\xi, \, (\xi\approx 0)$, whenever $x$
is infinitesimal). Therefore, by Theorem \ref{theorem1}, $f_1$ is
mu-differentiable.

In this case, the partial derivatives of $f$ are given by
$$\begin{array}{ll}
\displaystyle \left. \frac{\partial f}{\partial x_1} \right|_{x} &
\approx \displaystyle
\frac{f(x_1+\delta,x_2)-f(x_1,x_2)}{\delta}\\
                                                   & =
\displaystyle \frac{\sin(\epsilon x_1 + \epsilon
\delta)-\sin(\epsilon x_1)}{\epsilon\delta}\\[0.2cm]
                                                   & =
\displaystyle \frac{\cos(x^*)\cdot
\epsilon\delta}{\epsilon\delta}\\
                                                   &= \cos(x^*)\\
\end{array}$$
for some $x^* \in [\epsilon x_1, \epsilon x_1 + \epsilon \delta]$.
Therefore $x^*\approx 0$ and
$$\left. \frac{\partial f}{\partial x_1} \right|_{x} \approx 1.$$
In a similar way we might prove that
$$\left. \frac{\partial f}{\partial x_2} \right|_{x} =1.$$

\begin{definition} Let $f:{^*\mathbb R}^n\to {^*\mathbb R}$ be a
function and $a \in \mathbb R^n$ a vector. We say that $a$ is a
\emph{local m-minimizer} of $f$ if
$$f(x) \maior f(a) \, \mbox{ for all } \, x \in {^*B_r(a)}:= \{x
\in {^*\mathbb R}^n \, | \,\, |x-a|<r  \},$$
where $r\in \mathbb R$ is a positive real number. Analogously, we
define \emph{local m-maximizer} of $f$.
\end{definition}

Similarly to Lemma~\ref{ponte}, there exists a correspondence
between m-minimizers of internal functions and minimizers of
standard functions.

\begin{lemma}\label{ponte2}  If $f:{^*\mathbb R}^n\to {^*\mathbb
R}$ is mu-differentiable, then,
$$a \mbox{ is a m-minimizer of } f \mbox{ if and only if } a
\mbox{ is a minimizer of } st(f).$$
\end{lemma}

\begin{theorem} If $f:{^*\mathbb R}^n\to {^*\mathbb R}$ is a
mu-differentiable function and $a$ is a m-minimizer of $f$, then
$$\left. \frac{\partial f}{\partial x_i} \right|_{a}\approx 0 ,
\mbox{ for every }\, i=1,\ldots,n.$$
\end{theorem}
\begin{proof} If $a$ is a m-minimizer of $f$, by
Lemma~\ref{ponte2}, $a$ is a minimizer of $st(f)$. Therefore, for
all $i=1,\ldots,n$, $\displaystyle \left. \frac{\partial
st(f)}{\partial x_i} \right|_{a}=0$. By (\ref{eq3}), $$\left.
\frac{\partial f}{\partial x_i} \right|_{a}\approx 0.$$
\end{proof}

For example, the function $f$ defined in (\ref{eq2}) has no
m-minimums nor m-maximums.


\section*{Acknowledgment}

This work was supported by {\it Centre for Research on
Optimization and Control} (CEOC) from the ``Funda\c{c}\~{a}o para
a Ci\^{e}ncia e a Tecnologia'' (FCT), cofinanced by the European
Community Fund FEDER/POCI 2010.



\end{document}